\documentstyle[12pt, oldlfont]{article} 

\input{amssym.def}
\input{amssym.tex}

\begin{document}
\def\Ubf#1{{\baselineskip=0pt\vtop{\hbox{$#1$}\hbox{$\sim$}}}{}}
\def\ubf#1{{\baselineskip=0pt\vtop{\hbox{$#1$}\hbox{$\scriptscriptstyle\sim$}}}{}}
\def\R{{\Bbb R}}
\def\V{{\Bbb V}}
\def\N{{\Bbb N}}
\def\Q{{\Bbb Q}}

\title{Vaught's conjecture on analytic sets}         
\author{Greg Hjorth \footnote{Research partially supported by NSF grant DMS 96-22977}}        
\date{\today}          
\maketitle

{\bf $\S$0 Prehistory}

In rough historical these are the groups for which we know the topological 
Vaught conjecture:  

\medskip 

{\bf 0.1 Theorem} (Folklore) All locally compact Polish groups satisfy Vaught's 
conjecture -- that is to say, if $G$ is a locally compact Polish group acting 
continuously on a Polish space $X$ then either $|X/G|\leq\aleph_0$ or there is a 
perfect set of points with different orbits (and hence $|X/G|\geq 2^{\aleph_0}$). 

\medskip 

{\bf 0.2 Theorem} (Sami) Abelian Polish groups satisfy Vaught's conjecture. 

\medskip 

{\bf 0.3 Theorem} (Hjorth-Solecki) Invariantly metrizable and nilpotent Polish 
groups 
satisfy Vaught's conjecture. 

\medskip 

{\bf 0.4 Theorem} (Becker) Complete left invariant metric and solvable Polish 
groups satisfy Vaught's conjecture. 

\medskip 

In each of these case the result was shortly or immediately after extended to 
analytic sets. For this purpose let us write TVC$(G,\Ubf{\Sigma}^1_1)$ if whenever 
$G$ acts continuously on a Polish space $X$ and $A\subset X$ is $\Ubf{\Sigma}^1_1$ 
(or {\it analytic}) then either $|A/G|\leq\aleph_0$ or there is a 
perfect set of orbit inequivalent points in $A$. Thus we have TVC$(G,\Ubf{\Sigma}^1_1)$ 
for each of the group in the class mentioned in 0.1-0.4 above. 

On the other hand, and in contrast to the usual topological Vaught conjecture, that 
merely asserts that 0.1-0.4 hold for arbitrary Polish groups, it is known that 
TVC$(S_{\infty},\Ubf{\Sigma}^1_1)$ {\it fails}. 

Here it is shown that the presence of $S_{\infty}$ is a necessary condition for 
TVC$(G,\Ubf{\Sigma}^1_1)$ to fail: 

\medskip 

{\bf 0.5 Theorem} If $G$ is a Polish group on which the Vaught conjecture fails on 
analytic sets then there is a closed subgroup of $G$ that has $S_{\infty}$ as a continuous homomorphic image. 

\medskip 

The converse of 0.5 is known and by now considered trivial in light of 2.3.5 of 
\cite{beckerkechris}. Thus we have an exact characterization of 
TVC$(G,\Ubf{\Sigma}^1_1)$. If as widely believed the Vaught conjecture 
should fail for $S_{\infty}$ then this would as well characterize the groups for 
which the topological Vaught conjecture holds.

\newpage 

{\bf $\S$1 Preliminaries} 

\medskip 

All of this can be found in \cite{hjorthorbit}. 

\medskip

{\bf 1.1 Theorem} (Effros) Let $G$ be a Polish group acting continuously on 
a Polish space $X$ (in other words, let $X$ be a {\it Polish $G$-space}. 
For $x\in X$ we have $[x]_G\in \Ubf{\Pi}^0_2$ if and only if 
\[G\rightarrow [x]_G,\] 
\[g\mapsto g\cdot x\] 
is open. 

\medskip 

{\bf 1.2 Corollary} Let $G$ be a Polish group and  
$X$ a Polish $G$-space.  
Suppose that $[x]_G$ is $\Ubf{\Pi}^0_2$. 

Then for all $V$ containing the identity we may find open 
$U$ 
such that for all $ x'\in U\cap[x]_G$ and $U'\subset X$ open 
\[[x]_G\cap U'\cap U\neq\emptyset\]
implies that there exists $g\in V$ 
such that 
\[g\cdot x'\in U'.\] 

\medskip 

{\bf 1.3 Definition} Let $X$ be a Polish space and ${\cal B}$ a basis. Let ${\cal L}({\cal B})$ be the 
propositional language formed from the atomic propositions $\dot{x}\in U$, for $U\in{\cal B}$. 
Let ${\cal L}_{\infty,0}({\cal B})$ be the infinitary version, obtained by closing under negation and 
arbitrary disjunction and conjunction. 
$F\subset$ ${\cal L}_{\infty, 0}({\cal B})$ is a {\it fragment} if it is closed under subformulas and 
the finitary Boolean operations of negation and finite disjunction and finite conjunction. 

For a point $x\in X$ and $\varphi\in$${\cal L}_{\infty 0}({\cal B})$, 
we can then define 
$x\models \varphi$ by induction in the usual fashion, starting with 
\[x\models \dot{x}\in U\] 
if in fact $x\in U$. In the case that $X$ is a 
Polish $G$-space and $V\subset G$ open we may also define 
$\varphi^{\Delta V}$ by induction on the logical complexity of $\varphi$ so that 
in any generic extension in which $\varphi$ is hereditarily countable 
\[x\models \varphi^{\Delta V}\] 
if and only if 
\[\exists ^* g\in V (g\cdot x\models\varphi)\]
(where $\exists^*$ is the categoricity quantifier 
``there exists non-meagerly many'').

\medskip 

{\bf 1.4 Lemma} Let $X$ be a Polish $G$-space. ${\Bbb P}$ a forcing notion, $p\in {\Bbb P}$ 
a condition, $\sigma$ a ${\Bbb P}$-term. 
Suppose that ${\cal B}$ is a countable basis for $X$ and ${\cal B}_0$ a countable 
basis for $G$. Suppose that $G_0$ is a countable dense subgroup of $G$ and ${\cal B}$ is 
closed under $G_0$ translation and that ${\cal B}_0$ is closed under left and right 
$G_0$ translation. 
Suppose that 
\[p\Vdash_{\Bbb P}\sigma[\dot{G}]\in X\]
and that $p$ decides the equivalence class of $\sigma$ in the sense that 
\[(p,p)\Vdash_{{\Bbb P}\times{\Bbb P}}\sigma[\dot{G}_l]E_G\sigma[\dot{G}_r].\] 
Then there is a formula $\varphi_0$ and a fragment $F_0$ containing $\varphi_0$ 
so that: 

\leftskip 0.5in 

\noindent (i) $\{\{x\in X: x\models \psi^{\Delta V}\}: 
\psi\in F_0, V\in {\cal B}_0\}$ provides the basis for a topology $\tau_0(F_0)$, 
and in any generic extension in which $F_0$ becomes countable $(X,\tau_0(F_0))$ 
is a Polish $G$-space; 

\noindent (ii) $\varphi_0$ describes the equvalence class indicated by the triple 
$({\Bbb P},p,\sigma)$, in the sense that 
\[p\Vdash_{\Bbb P}\forall x\in X(x E_G\sigma[\dot{G}]\Leftrightarrow x\models
\varphi_0);\] 

\noindent (iii) and (ii) perserveres through all further forcing extensions, in that 
if $H\subset {\Bbb P}$ is $V$-generic below $p$ and $x=\sigma[H]$, then for all forcing notions ${\Bbb P}'\in V[H]$ 
\[ V[H]\models {\Bbb P}'\Vdash \forall y\in X(y E_G x\Leftrightarrow y\models \varphi_0).\] 

\leftskip 0in 

\medskip 

{\bf 1.5 Lemma} Let $G$ be a Polish group, $X$ a Polish $G$-space, $A\subset X$ a 
$\Ubf{\Sigma}^1_1$ set displaying a counterexample to TVC$(G,\Ubf{\Sigma}^1_1)$ -- 
so that $A/G$ has uncountably many orbits, but no perfect set of $E_G$-inequivalent 
points. 

Then for each ordinal $\delta$ there is a sequence $({\Bbb P}_{\alpha},
p_{\alpha},\sigma_{\alpha})_{\alpha <\delta}$ so that: 

\leftskip 0.5in 

\noindent (i) for each $\alpha<\delta$ 
\[(p_{\alpha},p_{\alpha})\Vdash_{{\Bbb P}_{\alpha}\times{\Bbb P}_{\alpha}}\sigma_{\alpha}[\dot{G}_l]E_G\sigma_{\alpha} [\dot{G}_r];\]

\noindent (ii) for each $\alpha<\beta< \delta$ 
\[(p_{\alpha},p_{\beta})\Vdash_{{\Bbb P}_{\alpha}\times{\Bbb P}_{\beta}}\neg(\sigma_{\alpha}[\dot{G}_l] E_G\sigma_{\beta} [\dot{G}_r]).\]

\leftskip 0in 

\newpage 

{\bf $\S$2 Proof} 

\medskip 

{\bf 2.1 Definition} $U$ is a {\it regular open} set if 
\[(\overline{U})^o=U\]
-- $U$ equals the interior of its closure. 
For $A$ a set let $RO(A)=(\overline{A})^o$. 

\medskip 

Note then that $RO(A)$ is always a regular open set. 

\medskip 

{\bf 2.2 Lemma} Let $G$ be a Polish group. For $V_0, V_1\subset G$ 
regular open sets, 
\[\{g\in G: V_0\cdot g= V_1\}\]
is a closed subset of $G$. \hfill ($\Box$) 

\medskip 

I need that the reader is willing to allow that we may speak of an 
$\omega$-model of set theory containing a Polish space, group, action, 
Borel set, and so on, provided suitable codes exist in the well founded 
part. Illfounded $\omega$-models are essential to the arguments below. 

In what follows let ZFC$^*$ be some large fragment of ZFC, at the 
very least 
strong enough to prove all the lemmas of $\S$1, but weak enough to admit a finite 
axiomatization.  

\medskip 

{\bf 2.3 Lemma} Let $M$ be an $\omega$-model of  ZFC$^*$. Let $X$, $G$, $G_0$, 
${\Bbb P}$, 
$F_0$, and so on, be as in 1.4 inside $M$. Suppose 
\[\pi:M\cong M\]
is an automorphism of $M$ fixing $X$, $G$, $G_0$, ${\Bbb P}$, $F_0$, $\varphi_0$, and 
all elements of ${\cal B}$ and ${\cal B}_0$. 
Suppose $H\subset$ Coll($\omega, F_0$) is $M$-generic and 
$x\in X^{M[H]}$ with 
\[x\models \varphi_0.\]  

Then there exists $\bar{g}\in G$ so that for all 
$\psi \in F_0$ and $V\in{\cal B}_0$ 
\[RO(\{g\in G_0: M[H]\models(g\cdot x\models \psi^{\Delta V})\})\bar{g}^{-1}= 
RO(\{g\in G_0: M[H]\models(g\cdot x\models \pi(\psi)^{\Delta V})\}).\] 

Proof. It suffices to find $g_0, g_1\in G$ so that 
\[RO(\{g\in G_0: M[H]\models(g\cdot x\models \psi^{\Delta V})\}){g}_0^{-1}= 
RO(\{g\in G_0: M[H]\models(g\cdot x\models \pi(\psi)^{\Delta V})\})g^{-1}_1\] 
for all $\psi$ and $V$. 

Let ${{\Bbb P}_0} $ be the forcing notion Coll$(\omega, F_0)$. 
Fixing $d_G$ a complete metric on $G$ we also build $h_i, h_i'\in G_0$, 
$\psi_i, \psi_i'\in F_0$, $W_i, W_i'\in {\cal B}_0$ so that 

\leftskip 0.5in 

\noindent (i) each $W_i$ is an open neighbourhood of the 
identity, $W_{i+1}\subset W_i$, $d_G(W_i)<2^{-i}$; 

\noindent (ii) $\pi(\psi_i)=\psi_i'$; 

\noindent (iii) $h_{2i}=h_{2i+1}$; $\forall g\in 
W_{2i+1}h_{2i}(d_G(g, h_{2i})<2^{-i})$; 

\noindent (iv)  $h_{2i+1}'=h_{2i+2}'$; $\forall g\in 
W_{2i+2}h_{2i+1}'(d_G(g, h_{2i+1}')<2^{-i})$; 

\noindent (v) $h_{i+1}\in W_ih_i$; $h_{i+1}'\in W_ih_i'$; 

\noindent (vi) $M[H]\models (h_i\cdot x\models (\psi_i)^{\Delta V_i})$; 

\noindent (vii) $M[H]\models (h_i'\cdot x\models (\psi_i')^{\Delta V_i})$; 

\noindent (viii) $M^{{\Bbb P}_0} $ satisfies that for all $y_0, y_1\in X$ 
all $\psi\in F_0$, and all $V\in {\cal B}_0$, 
if 
\[y_0\models \varphi_0\wedge (\psi_i)^{\Delta V_i}\]
\[y_1\models \varphi_0\wedge (\psi_i)^{\Delta V_i}\wedge \psi^{\Delta V}\] 
then 
\[y_0\models ((\psi_i)^{\Delta V_i}\wedge \psi^{\Delta V})^{\Delta W_i};\]

\noindent (ix) conversely $M^{{\Bbb P}_0} $ satisfies that for all $y_0, y_1\in X$, 
$\psi\in F_0$, $V\in {\cal B}_0$, 
if 
\[y_0\models \varphi_0\wedge (\psi_i')^{\Delta V_i}\]
\[y_1\models \varphi_0\wedge (\psi_i')^{\Delta V_i}\wedge \psi^{\Delta V}\] 
then 
\[y_0\models ((\psi_i')^{\Delta V_i}\wedge \psi^{\Delta V})^{\Delta W_i}.\] 

\leftskip 0in 

\noindent Note that (ix) actually follows from (viii), (ii), and the elementarity of 
$\pi$. 

Before verifying that we may produce $h_i, h_i'\in G_0$, 
$\psi_i, \psi_i'\in F_0$, $W_i, V_i\in {\cal B}_0$ as above, let us imagine that 
it is already completed and see how to finish. Using (iii) and (iv) we may obtain 
$g_0=$ lim$h_i$ and $g_1=$ lim$h_i'$. It suffices to check that for all 
\[ g\in RO(\{h\in G_0: M[H]\models(x\models \psi^{\Delta V})\}){g}_0^{-1}\]
we have 
\[g\in \overline{(\{h\in G_0: 
M[H]\models(h\cdot x\models \pi(\psi)^{\Delta V})\})} g_1^{-1}\]
(since the converse implication will be exactly symmetric). 

Then for sufficiently large $i$ we may choose a sufficiently small 
open neighbourhood 
$W$ of the identity and $\hat{g}\in G_0$ sufficiently close to $g$ 
so that $W\hat{g} W_i$ is an arbitrarily 
small neighbourhood of $g$ and 
\[M[H]\models (\hat{g} h_i\cdot x\models \psi^{\Delta V})\]
\[\therefore M[H]\models ( h_i\cdot x\models 
(\psi^{\Delta V})^{\Delta W\hat{g}})\] 
hence, as witnessed by $y=h_i\cdot x$ 
\[M^{{\Bbb P}_0}\models \exists y( y\models \varphi_0\wedge (\psi_i)^{\Delta V_i} 
\wedge (\psi^{\Delta V})^{\Delta W\hat{g}}),\]
\[\therefore M^{{\Bbb P}_0}\models \exists y( y\models \varphi_0\wedge (\psi_i')^{\Delta V_i} 
\wedge (\pi(\psi)^{\Delta V})^{\Delta W\hat{g}}),\] 
by elementarity of $\pi$, 
\[\therefore M[H]\models (h_i'\cdot x\models (\pi(\psi)^{\Delta V_i})^{\Delta W\hat{g}})
^{\Delta W_i})\]
by (ix), and so there exists some $\bar{g}\in W\hat{g} W_i$ so that 
\[M[H]\models (\bar{g}h_i'\cdot x\models \pi(\psi)^{\Delta V}).\] 
By letting $d_G(W\hat{g} W_i)\rightarrow 0$ and $h_i'\rightarrow g_1$ we get 
\[g\in \overline{\{h\in G_0: 
M[H]\models(x\models \pi(\psi)^{\Delta V})\}} g_1^{-1},\]
as required. 

We are left to hammer out the sequence. 

Suppose that we have $\psi_j, \psi'_j, W_j,V_j, h_j, h_j'$  for $j\leq 2i$. Immediately 
we may find $W_{2i+1}\subset W_{2i}$ giving (iii), and then by 1.2 and 1.4(i) we 
can produce $\psi_{2i+1}$, $V_{2i+1}$ satisfying (viii) and such that 
\[M[H]\models h_{2i}\cdot x=_{df} h_{2i+1}\cdot x\models (\psi_{2i+1})^{V_{2i+1}}.\] 
Then by considering that $\pi$ is elementary 
\[M^{{\Bbb P}_0} \models \exists y(y\models \varphi_0\wedge \pi(\psi_{2i})^{\Delta V_{2i}}
\wedge \pi(\psi_{2i+1})^{\Delta V_{2i+1}}).\]
Thus by (ix) we may find $h'\in G_0\cap W_{2i}$ so that 
\[M[H]\models (h'h_{2i}\cdot x\models \pi(\psi_{2i})^{\Delta V_{2i}}
\wedge \pi(\psi_{2i+1})^{\Delta V_{2i+1}}).\] 
In other words, by (ii), if we let $\psi_{2i+1}'=\pi(\psi_{2i+1})$ then 
\[M[H]\models (h'h_{2i}\cdot x\models (\psi_{2i}')^{\Delta V_{2i}}
\wedge (\psi_{2i+1}')^{\Delta V_{2i+1}}).\] 
Taking $h_{2i+1}'=h'h_{2i}'$ we complete the transition from $2i$ to $2i+1$. 

The  further step of producing $\psi_{2i+2}, \psi_{2i+2}', 
W_{2j+2}, h_{2j+2}$, $V_{2j+2}$ and $h_{2j+2}'$ is completely symmetrical. \hfill $\Box$ 

\medskip 

{\bf 2.4 Definition} $S_{\infty}$ {\it divides} a Polish 
group $G$ if there is a closed subgroup $H<G$ and a continuous onto homomorphism 
\[\pi:H\twoheadrightarrow S_{\infty}.\]
(By Pettis' lemma, any Borel homomorphism between Polish groups must be continuous.) 

\medskip 

{\bf 2.5 Lemma} $S_{\infty}$ divides Aut$({{\Bbb Q}} , <)$, the automorphism 
group of the rationals equipped with the usual linear ordering. 

\hfill ($\Box$) 

\medskip 

{\bf 2.6 Definition} For $X$, $G$, $F_0$, and so on, as in 1.4, ${\Bbb P}_0$= Coll$(\omega, 
F_0)$, $\psi_0, \psi_1\in F_0$, $V_0, V_1\in {\cal B}_0$, set 
\[(\psi_0, V_0)R (\psi_1, V_1)\]
if in $V^{{\Bbb P}_0}$ for all $x\models \varphi_0$ 
\[RO(\{g\in G_0: g\cdot x\models (\psi_0)^{\Delta V_0}\})\cap 
RO(\{g\in G_0: g\cdot x\models (\psi_1)^{\Delta V_1}\}\neq\emptyset).\] 
For $V\in {\cal B}_0$ let ${\cal B}(V)$ be the set of pairs $(\varphi, W)$ such that 
for all $\psi\in F_0$ and $W'\in {\cal B}_0$ 
\[V^{{\Bbb P}_0}\models \forall x_0\models \varphi_0\wedge \varphi^{\Delta W}
((\exists x_1\models \varphi_0\wedge \varphi^{\Delta W}\wedge \psi^{\Delta W'})
\Rightarrow x_0\models (\varphi^{\Delta W}\wedge \psi^{\Delta W'})^{\Delta V}).\]
In other words, ${\cal B}(V)$ corresponds to the basic open sets witnessing 1.2 
for $V$ in the topology $\tau_0(F_0)$. 

\medskip 

The next lemma states that if the equivalence class corresponding to $\varphi_0$ 
requires large forcing to be introduced then the formulas $\{\psi^{\Delta V}: 
\psi\in F_0, V\in{\cal B}_0\}$ have large $R$-discrete sets. 

\medskip 

{\bf 2.7 Lemma} Let $X$, $G$, $F_0$, ${\Bbb P}$, $\varphi_0$, 
and so on, be as in 1.4. Let $R$ be as in 2.6. Let $\kappa$ be a cardinal. 
Suppose no forcing notion of size less than $\kappa$ introduces a point in 
$X$ satisfying $\varphi_0$. 

Then there is no infinite $\delta<\kappa$ such that each ${\cal B}(V)$ for $V\in 
{\cal B}_0$ has a maximal $R$-discrete set of size $\leq \delta$. 

Proof. Suppose otherwise and choose large $\theta>\kappa$ so that $V_{\theta}\models$ 
ZFC$^*$ and choose an elementary substructure 
\[A {\prec} V_{\theta}\]
so that 
\[|A|=\delta,\]
\[\delta+1\subset A,\]
and $X$, $G$, $F_0$, $\varphi_0$, and so on, in $A$. Let $N$ be the transitive 
collapse of $A$ and 
\[\pi: N\rightarrow V_{\theta}\]
the inverse of the collapsing map. Set $\hat{{\Bbb P}}=\pi^{-1}({\Bbb P}_0)$ 
(where ${{\Bbb P}_0}$= Coll$(\omega, F_0)$), $ \hat{\varphi}_0=\pi^{-1}(\varphi_0)$, 
$\hat{F}_0=\pi^{-1}(F_0)$, 
choose 
\[\hat{H}\subset \hat{\Bbb P},\]
\[H\subset {\Bbb P}_0\]
to be $V$-generic, and choose $\hat{x}\in N[\hat{H}]$ and $x\in V[H]$ so that 
\[N[\hat{H}]\models (\hat{x}\models \hat{\varphi}_0),\]
\[V[H]\models (x\models\varphi_0).\] 
It suffices to show 
\[\hat{x}E_Gx.\]
 
As in the proof of 2.3 find $h_i, h_i'\in G_0$, $\psi_i\in F_0$, 
$\psi_i'\in {\hat{F}}_0$, $V_i, V'_i \in{\cal B}_0$, $W_i\in {\cal B}_0$ 
and $U_i\subset X$ basic open so that : 

\leftskip 0.5in 

\noindent (i) $W_{i+1}\subset W_i$, $W_i=(W_i)^{-1}$, $d_G(W_i)<2^{-i},$ $1_G\in W_i$; 
$U_{i+1}\subset U_i$,  
$d_X(U_i)<2^{-i}$; 

\noindent (ii) $\pi(\psi_i')=\psi_i$; 

\noindent (iii) $\forall g\in (W_{2i+1})^3 h_{2i}(d_G(g, h_{2i})<2^{-i})$; 
$h_{2i+1}=h_{2i}$; 

\noindent (iv) $\forall g\in (W_{2i+2})^3 h_{2i+1}'(d_G(g, h_{2i+1})<2^{-i})$; 
$h_{2i+2}'=h_{2i+1}'$; 

\noindent (v) $h_{i+1}\in (W_i)^3 h_i$, $h_{i+1}'\in (W_i)^3 h_i'$; 

\noindent (vi) $V[H]\models (h_i\cdot {x}\models (\psi_i)^{\Delta V_i})$; 

\noindent (vii) $N[\hat{H}]\models (h_i'\cdot \hat{x}\models (\psi_i')^{\Delta V_i'})$;

\noindent (viii) $V^{{\Bbb P}_0}\models (\psi_i, V_i)\in {\cal B}(W_i)$; 

\noindent (ix) $N^{{\Bbb P}_0}\models (\psi_i', V_i')\in {\cal B}(W_i)$; 

\noindent (x) $(\pi(\psi_i'), V_i') R(\psi_i, V_i)$; 

\noindent (xi) $h_i\cdot x, h_i'\cdot \hat{x} \in U_i$. 

\leftskip 0in 

Granting all this may be found we finish quickly. By (iii) and (iv) we 
get $g_0$= lim$h_i$ and  $g_1$= lim$h_i'$, whence 
\[g_0\cdot x=g_1\cdot \hat{x}\]
by (xi). This would contradict ${\hat{\Bbb P}}$ being too small to introduce 
a representative of $[x]_G$. 

So instead suppose we have built $V_j, V_j', \psi_j$ and so on for $j\leq 2i$ and 
concentrate on trying to show that we may continue the construction up to 
$2i+2$. 

First choose $W_{2i+1}\subset W_{2i}$ in accordance with (i) and (iii) and then for 
(xi) and (i) choose $U_{2i+1}\subset U_{2i}$ containing $h_{2i}\cdot x (=_{df} 
h_{2i+1}\cdot x)$ with $d_X(U_{2i+1})<2^{-2i-1}$. Then by 1.2 we may choose 
$V_{2i+1}, \psi_{2i+1}$ as indicated at (vi) and (viii). 

On the $N$ side we use the assumption on $R$ to find $V_{2i+1}'$ and $\psi_{2i+1}'$ in 
$N$ so that 
\[N^{\hat{\Bbb P}_0}\models (\psi_{2i+1}', V_{2i+1}')\in {\cal B}(W_{2i+1})\]
and 
\[(\pi(\psi_{2i+1}'), V_{2i+1}') R (\psi_{2i+1}, V_{2i+1}).\] 
Unwinding the definitions gives 
\[V^{{\Bbb P}_0}\models (y\models \varphi_0 \wedge \pi(\psi_{2i}')^{\Delta V_{2i}'})
\Rightarrow y\models ((\psi_{2i})^{\Delta V_{2i}}
\wedge \pi(\psi_{2i}')^{\Delta V_{2i}'})^{\Delta W_{2i}},\] 
\[V^{{\Bbb P}_0}\models (y\models \varphi_0 \wedge (\psi_{2i})^{\Delta V_{2i}})
\Rightarrow y\models ((\psi_{2i})^{\Delta V_{2i}}
\wedge (\psi_{2i+1})^{\Delta V_{2i+1}})^{\Delta W_{2i}},\] 
\[V^{{\Bbb P}_0}\models (y\models \varphi_0 \wedge (\psi_{2i+1})^{\Delta V_{2i+1}})
\Rightarrow y\models ((\psi_{2i+1})^{\Delta V_{2i+1}}
\wedge \pi(\psi_{2i+1}')^{\Delta V_{2i+1}'})^{\Delta W_{2i+1}}.\] 
In particular, assuming without loss of generality that 
\[(\psi_{2i+1})^{\Delta V_{2i+1}}\Rightarrow \dot{x}\in U_{2i+1}\]
we have 
\[V^{{\Bbb P}_0}\models (y\models \varphi_0 \wedge \pi(\psi_{2i}')^{\Delta V_{2i}'})
\Rightarrow y\models ((\psi_{2i+1}')^{\Delta V_{2i+1}{'}}
\wedge \dot{x}\in U_{2i+1})^{\Delta (W_{2i})^3}.\] 
Thus by elementarity of $\pi$ we may find $h'\in (W_{2i})^3\cap G_0$ so that $h'h_{2i}'\cdot \hat{x}\in 
U_{2i+1}$ and 
\[N[\hat{H}]\models (h'h_{2i}'\cdot \hat{x}\models (\psi_{2i+1}')^{\Delta V_{2i+1}'}).\] 
Then setting $h_{2i+1}=h'h_{2i}$ completes the transition from $2i$ to $2i+1$. 

The step from $2i+1$ to $2i+2$ is similar. \hfill $\Box$ 

\medskip 

We need a fact from infinitary model theory. 

\medskip 

{\bf 2.8 Theorem} Let $\varphi\in {\cal L}_{\omega_1,\omega}$ and suppose 
\[N\models \varphi\]
and $P$ is a predicate in the language of $N$ with 
\[|(P)^N|\geq \beth_{\aleph_1}.\] 
Then $\varphi$ has a model with generating indiscernibles in $P$. 

More precisely there is a model $M$ with language ${\cal L}^*\supset 
{\cal L}(N)$, ${\cal L}^*$ having a new symbol $<$, along with new 
function symbols of the form $f_{\hat{\varphi}}$ for $\hat{\varphi}$ in the 
fragment of ${\cal L}(N)_{\omega_1,\omega}$ generated by $\varphi$, 
and distinguished elements $(c_i)_{i\in \N}$, 
so that: 

\leftskip 0.5in 

\noindent (i) $(<)^M$ linearly orders $(P)^M$; 

\noindent (ii) each $f_{\hat{\varphi}}$ is a Skolem function for 
$\hat{\varphi}$; 

\noindent (iii) $M$ is the Skolem hull of $\{c_i:i\in\N\}$ (under the 
functions of the form  $f_{\hat{\varphi}}$); 

\noindent (iv) each $c_i\in (P)^M$; 

\noindent (v) for all $\psi$ in the fragment of ${\cal L}^*_{\omega_1,\omega}$ 
generated by $\varphi$ and $i_1<i_2<...<i_n$, $j_1<...<j_n$ in $\N$ 
\[M\models \psi(c_{i_1},c_{i_2},...,c_{i_n})\Leftrightarrow 
\psi(c_{j_1},c_{j_2},...,c_{j_n});\] 

\noindent (vi) $M\models \varphi$. 

\leftskip 0in 

See \cite{keisler}. \hfill ($\Box$) 

\medskip 

{\bf 2.9 Theorem} Let $G$ be a Polish group for which TVC($G,\Ubf{\Sigma}^1_1$) fails. 
Then $S_{\infty}$ divides $G$. 

Proof. Choose some Polish $G$-space $X$ witnessing the failure of TVC($G,\Ubf{\Sigma}^1_1$). Following 1.5 we may find some $({\Bbb P}, p,\sigma)$ 
introducing an equivalence class as in 1.4 that may not be produced by a forcing 
notion of size less than $\beth_{\aleph_1}$. Fix $\varphi_0$, ${\cal B}$, 
${\cal B}_0$, $F_0$, $G_0$, and so on, as in 1.4, so that in all generic extensions $V[H]$ of $V$ 
\[V[H]\models p\Vdash_{\Bbb P}\forall y\in X(yE_G\sigma[\dot{G}]\Leftrightarrow 
y\models \varphi_0).\] 

Let $V_{\theta}$ be large enough to contain $X$, $G$, $\varphi_0$, and so on, and 
satisfy ZFC$^*$. By 2.7 choose $P\subset V_{\theta}$ to be of size $\beth_{\aleph_1}$ and 
$R$-discrete (or more precisely, 
so for all $(\psi, V)\neq(\psi',V')\in P$ we have for any $V$-generic $H\subset$ 
Coll$(\omega, F_0)$ that $V[H]\models \neg ((\psi, V)R(\psi',V'))$). 
Applying 2.8 to $N=(V_{\theta}; \in, P, X, G, G_0, \varphi_0,...)$ and 
we may obtain an $\omega$-model 
with indiscernibles $(\psi_q, V_q)_{q\in {\Bbb Q}}$ in $B^M$. Let $H\subset$ 
Coll($\omega, (F_0)^M$) be $M$-generic. Choose $x\in M[H]$ so that 
\[M[H]\models (x\models \varphi_0).\] 

All this granted we may define $G_1$ to be the set of $\bar{g}\in G$ so that for all $q\in 
{\Bbb Q}$ there exists $r\in{\Bbb Q}$ with 
\[RO(\{g\in G_0: M[H]\models(g\cdot x\models (\psi_q)^{\Delta V_q})\})\bar{g}^{-1}= 
RO(\{g\in G_0: M[H]\models(g\cdot x\models (\psi_r)^{\Delta V_r})\})\] and for  
$q\in 
{\Bbb Q}$ there exists $r\in{\Bbb Q}$ with 
\[RO(\{g\in G_0: M[H]\models(g\cdot x\models (\psi_q)^{\Delta V_q})\})\bar{g}= 
RO(\{g\in G_0: M[H]\models(g\cdot x\models (\psi_r)^{\Delta V_r})\}).\] 

$G_1$ is $\Ubf{\Pi}^0_2$ in $G$, by 2.2 and since $\bar{g}$ is in $G_1$ if and only if the 
following four conditions hold: 

\leftskip 0.5in 

\noindent (i) for all $q\in 
{\Bbb Q}$ there exists $r\in{\Bbb Q}$ with 
\[RO(\{g\in G_0: M[H]\models(g\cdot x\models (\psi_q)^{\Delta V_q})\})\bar{g}^{-1}\cap 
RO(\{g\in G_0: M[H]\models(g\cdot x\models (\psi_r)^{\Delta V_r})\})\neq\emptyset,\]

\noindent (ii) for all $q,r\in 
{\Bbb Q}$ 
\[RO(\{g\in G_0: M[H]\models(g\cdot x\models (\psi_q)^{\Delta V_q})\})\bar{g}^{-1}\cap 
RO(\{g\in G_0: M[H]\models(g\cdot x\models (\psi_r)^{\Delta V_r})\})\neq\emptyset\]
implies 
\[RO(\{g\in G_0: M[H]\models(g\cdot x\models (\psi_q)^{\Delta V_q})\})\bar{g}^{-1}= 
RO(\{g\in G_0: M[H]\models(g\cdot x\models (\psi_r)^{\Delta V_r})\});\]

\noindent (iii) for all $q\in 
{\Bbb Q}$ there exists $r\in{\Bbb Q}$ with 
\[RO(\{g\in G_0: M[H]\models(g\cdot x\models (\psi_q)^{\Delta V_q})\})\bar{g}\cap 
RO(\{g\in G_0: M[H]\models(g\cdot x\models (\psi_r)^{\Delta V_r})\})\neq\emptyset,\]

\noindent (iv) for all $q,r\in 
{\Bbb Q}$ 
\[RO(\{g\in G_0: M[H]\models(g\cdot x\models (\psi_q)^{\Delta V_q})\})\bar{g}\cap 
RO(\{g\in G_0: M[H]\models(g\cdot x\models (\psi_r)^{\Delta V_r})\})\neq\emptyset\]
implies 
\[RO(\{g\in G_0: M[H]\models(g\cdot x\models \psi_q^{\Delta V_q})\})\bar{g}^{-1}= 
RO(\{g\in G_0: M[H]\models(g\cdot x\models (\psi_r)^{\Delta V_r})\}).\] 

\leftskip 0in 

\noindent Since $G_1$ is a $\Ubf{\Pi}^0_2$ subgroup of $G$ it must be closed. 

For $g\in G_1$ we may define the permutation $\hat{\pi}(g)$ of ${\Bbb Q}$ by the 
specification that for all $q\in {\Bbb Q}$ 
\[(\hat{\pi}(g))(q)=r\] 
if and only if $r$ is as above in the definition of $G_1$. This is well defined by 
the $R$-discreteness of the set $(P)^M$. 

Now let $G_2$ be the set of $g\in G_1$ such that $\hat{\pi}(g)$ defines an 
automorphism of the structure $({\Bbb Q}, <)$. $G_2$ is a closed subgroup of 
$G_1$ and hence $G$. Since every order preserving permutation of the indiscernibles 
induces an automorphism of $M$ the map 
\[\hat{\pi}:G_2\rightarrow {\rm Aut}({\Bbb Q}, <)\]
is onto by 2.3. Then by 2.5 $S_{\infty}$ divides $G$. \hfill $\Box$ 

\medskip 

{\bf 2.10 Conjecture} Assume AD$^{L(\R)}$. Let $G$ be a Polish group, 
$X$ a Polish $G$-space, $A\subset X$ in $\Ubf{\Sigma}^1_1$, and suppose in 
${L(\R)}$ there is an injection 
\[i:A/G\hookrightarrow 2^{<\omega_1}.\]
Then there is a Polish $S_{\infty}$-space $Y$ and a $\Ubf{\Sigma}^1_1$ set 
$B\subset Y$ and a bijection 
\[\pi:A/G\cong B/S_{\infty}.\]

6363 MSB

Mathematics

UCLA

CA90095-1555

greg@math.ucla.edu

\end{document}